\documentclass{article}
\usepackage{amsmath,amssymb,amsthm,latexcad}

\title{\bf{Some properties of the value function and its level
sets for affine control systems with quadratic cost}}
\author{E. Tr\'elat \\ \\
Universit\'e de Bourgogne, Laboratoire de Topologie, \\
UMR 5584 du CNRS, BP47870, 21078 Dijon Cedex, France  \\
e-mail : trelat@topolog.u-bourgogne.fr}
\date{}
\vspace{1cm}

\newtheorem{thm}{Theorem}[section]
\newtheorem{cor}[thm]{Corollary}
\newtheorem{prop}[thm]{Proposition}
\newtheorem{lem}[thm]{Lemma}

\theoremstyle{remark}
\newtheorem{rem}{Remark}[section]

\theoremstyle{definition}
\newtheorem{defi}{Definition}[section]

\theoremstyle{definition}
\newtheorem{ex}{Example}[section]

\renewcommand{\geq}{\geqslant}
\renewcommand{\leq}{\leqslant}
\renewcommand{\l}{\left}
\renewcommand{\r}{\right}
\newcommand{\f}[2]{\frac{#1}{#2}}
\newcommand{\inv}[1]{\frac{1}{#1}}
\renewcommand{\it}{\textit}
\renewcommand{\rm}{\textrm}
\newcommand{\tto}{\longrightarrow}
\renewcommand{\to}{\rightarrow}

\def\N{\rm{I\kern-0.21emN}}

\def\R{\rm{I\kern-0.21emR}}
\def\H{\rm{I\kern-0.21emH}}
\def\K{\rm{I\kern-0.24emK}}
\def\P{\rm{I\kern-0.21emP}}
\def\Q{\rm{l\kern-0.5emQ}}

\newcommand{\no}{\noindent}
\newcommand{\ds}{\displaystyle}

\newcommand{\esp}{\,\,\,}

\begin{document}
\maketitle

\begin{abstract}
\no Let $T>0$ fixed.
We consider the optimal control problem for analytic
affine systems~:
$\ds{\dot{x}=f_0(x)+\sum_{i=1}^m u_if_i(x)}$, with a cost of the
form~: $\ds{C(u)=\int_0^T \sum_{i=1}^m u_i^2(t)dt}$.
For this kind of systems we prove
that if there are no minimizing
abnormal extremals then the value function $S$
is subanalytic. Secondly we prove that if there exists an
abnormal minimizer of corank 1 then the set of end-points of
minimizers at cost fixed is tangent to a given hyperplane.
We illustrate this situation in sub-Riemannian geometry.
\end{abstract}

\paragraph{Key words :} optimal control,
value function, abnormal minimizers,
subanalyticity, sub-Riemannian geometry.

\section{Introduction} 
Let $M$ be an analytic Riemannian $n$-dimensional manifold and
$x_0\in M$. Consider the following \it{control system}~:
\begin{equation} \label{syst0}
{\dot x}(t)=f(x(t),u(t)),\esp
x(0)=x_0  
\end{equation}
where $f=M\times \R^m \longrightarrow M$ is an analytic
function, and the set of controls $\Omega$ is
a subset of measurable mappings defined on
$[0,T(u)]$ and taking their values in $\R^m.$
The system is said to be \it{affine} if~:
\begin{equation}  \label{affine0}
f(x,u)=f_0(x)+\sum_{i=1}^m u_if_i(x)
\end{equation}
where the $f_i$'s are analytic vector fields on $M$.

Let $T>0$ fixed. We consider the \it{end-point mapping}
$E:u\in \Omega \longmapsto x(T,x_0,u),$ where $x(t,x_0,u)$
is the solution of (\ref{syst0}) associated to $u\in \Omega$ and
starting from $x_0$ at $t=0.$ We endow the set of controls defined
on $[0,T]$ with the $L^2$-norm topology. A trajectory
${\widetilde x} (t,x_0,{\widetilde u})$ denoted in short
${\widetilde x}$ is said to be \it{singular or abnormal} on
$[0,T]$ if ${\widetilde u}$ is a singular point of the end-point
mapping, i.e, the Fr\'echet derivative of $E$ is not surjective
at ${\widetilde u}$ ; otherwise it is said \it{regular}.
We denote by $Acc(T)$ the set of end-points
at $t=T$ of solutions of (\ref{syst0}), $u$ varying in $\Omega$.
The main problem of control theory is to study $E$ and $Acc(T)$.
Note that the latter is not bounded in general. In \cite{LM} one
can find sufficient conditions so that $Acc(T)$ is compact, or
has non-empty interior. Theorem \ref{compact} of this article
states such a result for affine systems.

Consider now the following \it{optimal control problem}~: among
all trajectories of (\ref{syst0}) steering $0$ to $x\in Acc(T)$,
find a trajectory minimizing the \it{cost function}~:
$\ds{C(u)=\int_0^T f^0(x_u(t),u(t))dt}$, where $f^0$ is analytic.
Such minimizers do not necessary exist ; the main
argument to prove existence theorems is the lower
semi-continuity of the cost function, see \cite{F} or \cite{LM}.
If $x\in Acc(T)$, we set $S(x)=\inf\{C(u)\esp /\esp E(u)=x\}$,
otherwise $S(x)=+\infty$ ; $S$ is called
the \it{value function}. In general $f^0$ is choosen
in such a way that
the value function has a physical meaning~: for instance
the action in classical mechanics or in optics, the
(sub)-Riemannian distance in (sub)-Riemannian geometry.
We are interested in the regularity of the value
function and the structure of its level
sets. In (sub)-Riemannian geometry level sets of the
distance are (sub)-Riemannian spheres. To describe these objects
we need a category of sets which are stable under set operations
and under proper analytic maps.

An important example of such a category is the one of
\it{subanalytic sets} (see \cite{Gab}). They have been utilized
by several authors in order to construct an optimal synthesis
or to describe $Acc(T)$ (see \cite{Bru}, \cite{Su}). Unfortunately
this class is not wide enough~: in
\cite{LoS}, the authors exhibit examples of control systems in
which neither $S$ nor $Acc(T)$ are subanalytic. However Agrachev
shows in \cite{Ag} (see also \cite{AgS}, \cite{J}) that if there
are no abnormal minimizers then the sub-Riemannian distance
is subanalytic in a pointed neighborhood of $0$, and hence
sub-Riemannian spheres of small radius are subanalytic.
Following his
ideas, we extend this result to affine control systems with
quadratic cost (Theorem \ref{subanalytic} and corollaries).

Abnormal minimizers are responsible for a phenomenon of
\it{non-properness} (Proposition \ref{nonproper}),
which geometrically implies the following property~:
under certain assumptions the level sets of the value function
are tangent to a given hyperplane at the end-point of the
abnormal minimizer (Theorem \ref{tangent}).
This result was first stated in \cite{BT} for
sub-Riemannian systems to illustrate the Martinet situation.

An essential reasoning we will use in the proofs of these results is
the following (see notably Lemma
\ref{AV0}). We shall consider sequences of minimizing controls
$u_n$ associated to \it{projectivized Lagrange multipliers}
$(p_n(T),p^0_n)$, so that we have (see section 2)~:
\begin{equation} \label{explic}
p_n(T)dE(u_n)=-p^0_nu_n
\end{equation}
Since $(u_n)$ is bounded in $L^2$ we shall assume that $u_n$
\it{converges weakly} to $u$. To pass to the limit in
(\ref{explic}) we shall prove some regularity properties of the
end-point mapping $E$ (section 3).
Contrarily to the sub-Riemannian case the
strong topology on $L^2$ is not adapted in general for affine
systems, whereas the weak topology gives nice compacity
properties of the set of minimizing controls (see Theorem
\ref{compactcontinu}).
\\

The outline of the paper is as follows~: in section 2 we recall
definitions of subanalytic sets and the Maximum Principle. In §3
we state some basic results on the regularity of the end-point
mapping. Section 4 is devoted to continuity and subanalyticity of
the value function $S$. Finally, in §5 the shape of the level
sets of the value function in presence of abnormal minimizers is
investigated. We illustrate this situation in sub-Riemannian
geometry.
\\

I would like to thank A. Agrachev for many advices and remarks
which helped me in this work.


\section{Preliminaries}
\subsection{Subanalytic sets}
Recall the following definitions, that can be found in
\cite{hardt}, \cite{hironaka}.

\begin{defi}
Let $M$ be a finite dimensional real analytic manifold.
A subset
$A$ of $M$ is called \it{semi-analytic} iff, for every $x$ in $M$,
we can find a neighborhood $U$ of $x$ in $M$ and $2pq$ real
analytic functions $g_{ij}$, $h_{ij}$
($1\leq i\leq p$ and $1\leq j\leq q$) such that
$$A\cap U=\overset{p}{\underset{i=1}{\bigcup}}
 \{ y\in U \esp / \esp g_{ij}(y)=0 \esp
\rm{and} \esp h_{ij}(y)>0 \rm{ for }j=1\ldots q \} $$
We let $SEM(M)$ denote the family of semi-analytic subsets of $M$.
\end{defi}

Unfortunately proper analytic images of semi-analytic sets are
not in general semi-analytic. Hence this class must be
extended~:

\begin{defi}
A subset
$A$ of $M$ is called \it{subanalytic} iff, for every $x$ in $M$,
we can find a neighborhood $U$ of $x$ in $M$ and $2p$ pairs
$(\phi_i^\delta, A_i^\delta)$ ($1\leq i\leq p$ and $\delta
=1,2$), where $A_i^\delta \in SEM(M_i^\delta)$ for some real
analytic manifolds $M_i^\delta$, and where the maps
$\phi_i^\delta : M_i^\delta \to M$ are proper analytic, such
that
$$A\cap U= \overset{p}{\underset{i=1}{\bigcup}}
( \phi_i^1 (A_i^1) \backslash \phi_i^2(A_i^2)) $$
We let $SUB(M)$ denote the family of subanalytic subsets of $M$.
\end{defi}

The class of subanalytic sets is closed under
union, intersection, complement, inverse image of analytic
maps, image of proper analytic maps. Moreover they are
\it{stratifiable}. Recall the following~:

\begin{defi}
Let $M$ be a differentiable manifold. A stratum in $M$ is a
locally closed submanifold of $M$.

A locally-finite partition ${\mathcal{S}}$ of $M$ is called a
stratification of $M$ if each $S$ in ${\cal{S}}$ is a stratum such
that~:
$$\forall T\in {\cal{S}} \quad T\cap \rm{Fr }S \neq \emptyset
\Rightarrow T \subset \rm{Fr }S \rm{ and } \rm{dim }T<\rm{dim
S}$$
\end{defi}

Endly, a map $f:M\to N$ between two manifolds is called
\it{subanalytic} if its graph is a subanalytic set of $M\times
N$.

The basic property of subanalytic functions which makes them
useful in optimal control theory is the following.
It can be found in \cite{Tamm}.

\begin{prop} \label{tamm137}
Let $M$ and $N$ denote finite dimensional real analytic
manifolds, and $A$ be a subset of $N$.
Given subanalytic maps $\phi : N \tto M$ and $f:N\tto \R$,
we define~:
$$\forall x\in M \quad \psi(x)=\inf \{f(y) \esp / \esp y\in \phi
^{-1}(x) \cap A \}$$
If $\phi_{/\bar{A}}$ is proper, then $\psi$ is subanalytic.
\end{prop}


\subsection{Maximum principle and extremals}
According to the \it{weak maximum principle} \cite{P} the minimizing
trajectories are among the singular trajectories of the end-point
mapping of the \it{extended system} in $M\times \R$~:
\begin{eqnarray}  \label{extended}
{\dot x} (t)& = & f(x(t),u(t))  \\ 
{\dot x}^0 (t) & = & f^0 (x(t),u(t))\nonumber
\end{eqnarray}

\noindent They are called \it{extremals}. If $E$ and $C$ are
differentiable, then there exists a
\it{Lagrange multiplier} $(p(T),p^0)$ (defined up to a scalar)
such that~:
\begin{equation} \label{Laggen}
p(T)dE(u)=-p^0dC(u)
\end{equation}
where $dE(u)$ (resp. $dC(u)$) denotes the differential of $E$
(resp. $C$) at $u$.
Moreover, $(x(T),p(T))$ is the end-point of the solution
of the following equations~:
\begin{eqnarray}  \label{constraint}
{\dot x}={\displaystyle {\partial H \over \partial p}}\ 
,\ \ {\dot p}=-{\displaystyle {\partial H \over \partial x}}
\ ,\ \ {\displaystyle {\partial H \over \partial u}}=0
\end{eqnarray}

\noindent where $H=<p,f(x,u)>+p^0 f^0 (x,u)$ is the
\it{Hamiltonian}, $p$ is the adjoint vector, $<,>$
the inner product on $M$ and $p^0$ is a constant.
The abnormal
trajectories correspond to the case $p^0=0$ and their role
in the optimal control problem has to be analyzed. The extremals
with $p^0\neq 0$ are said normal. In this case $p^0$ is usually
normalized to
$-\inv{2}$. We will use this normalization to prove Theorem
\ref{subanalytic}. To prove Theorem \ref{tangent} we will use
another normalization by considering \it{projectivized Lagrange
multipliers}, i.e. $(p(T),p^0)\in P(T^*M)$. We say that an
extremal has \it{corank 1} if it has an unique projectivized
Lagrange multiplier.

\paragraph{Affine systems}
Consider in particular analytic affine control systems on
$M$~:
\begin{equation} \label{aff0}
\dot{x}(t)=f_0(x)+\sum_{i=1}^m u_if_i(x),\esp
x(0)=0 
\end{equation}
where the $f_i$'s are analytic vector fields, with the problem
of minimizing the following cost~:
\begin{equation} \label{cost}
C(u)=\int_0^T \sum_{i=1}^m u_i^2(t)dt 
\end{equation}
The Hamiltonian is~:
$$H(x,p,u)=<p,f_0(x)+\sum_{i=1}^m u_if_i(x)>+p^0 \sum_{i=1}^m
u_i^2 $$

\paragraph{Parametrization of normal extremals}
We suppose $p^0=-\inv{2}$.
Then normal controls can be computed from
the equation~: ${\displaystyle {\partial H \over \partial u}}=0$,
and we get~:
\begin{equation} \label{norcont}
\forall i=1\ldots m \quad u_i=<p,f_i(x)>
\end{equation}
Putting in system (\ref{aff0}), we get an analytic differential
system in $T^*M$ parametrized by the initial condition $p(0)$.
From the general theory of ordinary differential equations we
know that solutions depend \it{analytically} on their initial
condition. Denote such a solution by
$(x_{p(0)},p_{p(0)})$. Let $u_{p(0)}=({<p_{p(0)},f_1(x_{p(0)})>,}
{\ldots,} <p_{p(0)},f_m(x_{p(0)})>)$~; from (\ref{norcont})
$u_{p(0)}$ is a \it{normal control} associated to $x_{p(0)}$. Now
we can define~:

\begin{defi} \label{fi}
The mapping $\Phi~:\begin{array}{rcl}
T^*_{x_0}M&\longrightarrow &L^2([0,T],\R^m) \\
p(0)&\longmapsto &u_{p(0)} \end{array}$ is \it{analytic}.
\end{defi}

\no This mapping will be useful to check subanalyticity of the
value function in section 4.


\section{Regularity of the end-point mapping}
Let  $M$ be an analytic complete
n-dimensional Riemannian manifold and $x_0 \in M$.
Our point of view is local and we can assume~: $M=\R^n, x_0=0$.
We only consider \it{analytic affine control systems}
(\ref{aff0}). The statements in this section except for
Proposition \ref{contdiff} are quite standard, and we include
proofs only for convenience of the reader.

\subsection{The end-point mapping}
Let $T>0$ and $x_u$ be the solution, if exists, of the controlled
system~:
\begin{equation*}
\dot{x}_u=f_0(x)+\sum_{i=1}^mu_if_i(x),\esp
x_u(0)=0
\end{equation*}
where $u=(u_1,\ldots,u_m)\in L^2([0,T],\R^m)$.
Since we allow discontinuous controls, the meaning of solution of
the previous differential system has to be clarified. In fact
this means that the following integral equation holds~:
$$\forall t\in [0,T]\quad x_u(t)=\int_0^Tf_0(x_u(\tau))+
\sum_{i=1}^mu_i(\tau)f_i(x_u(\tau)) \, d\tau$$

\begin{defi}
The \it{end-point mapping} is~:
$$\begin{array}{rcl}
E : \Omega & \tto & \R^n \\
u & \longmapsto & x_u(T)
\end{array}$$
where $\Omega\subset L^2([0,T],\R^m)$ is the domain of $E$, that
is the subset of controls $u$ such that $x_u$ is well-defined on
$[0,T]$.
\end{defi}

$E$ is not defined on the whole $L^2$ because of
\it{explosion phenomena}. For example consider the system
$\dot{x}=x^2+u$ ; then $x_u$ is not defined on $[0,T]$ for $u=1$
if $T\geq \f{\pi}{2}$. Anyway we have the following~:

\begin{prop} \label{existence}
Let $T>0$ fixed. We consider the analytic control system
(\ref{aff0}). Then the domain $\Omega$ of $E$ is open
in $L^2([0,T],\R^m)$.
\end{prop}

\begin{proof}[Proof of the Proposition]
It is enough to prove the following statement~:
\begin{quote}
If the trajectory $x_u$ associated to $u$ is well-defined on
$[0,T]$, then the same is true for any control in a neighborhood
of $u$ in $L^2([0,T],\R^m)$.
\end{quote}
Let $V$ be a bounded open subset of $\R^n$ such that~:
$\forall t \in [0,T] \quad x_u(t) \in V$. Let $\theta \in
C^\infty (\R^n , [0,1])$ with compact support $K$ such that $\theta
=1$ on $V$. We can assume that $K=\bar{B}(0,R)=\{ x\in \R^n \, /
\, ||x|| \leq R \}$. For $i=0 \ldots m$
we set $\tilde{f}_i=\theta f_i$.
Then it is clear that $x_u$ is also solution of~:
$\dot{x}=\tilde{f}_0(x)+\sum u_i \tilde{f}_i(x)$.
For all $v\in L^2$ let $\tilde{x}_{v}$ be the solution of
$\dot{\tilde{x}}_{v}=\tilde{f}_0(\tilde{x}_{v}) +
\displaystyle{\sum_{i=1}^m v_{i} \tilde{f}_i
(\tilde{x}_{v})}$, $\tilde{x}_{v}(0)=0$.
We will prove that $\tilde{x}_{v} =
x_{v}$ in a small enough neighborhood of $u$.

\begin{lem}
The $\tilde{f}_i$ are globally lipschitzian on $\R^n$, that is~:
$$\exists A>0 \esp / \esp \forall i\in \{0,\ldots,m\}
\esp \forall y,z \in \R^n \quad 
||\tilde{f}_i(y)-\tilde{f}_i(z)|| \leq A||y-z||$$
\end{lem}

\begin{proof}[Proof of the Lemma]
Let $i\in \{0,\ldots,m\}$. $\tilde{f}_i$ is $C^1$ hence is locally
lipschitzian at any point~:
$$\forall x\in \bar{B}(0,2R) \quad \exists \rho_x, A_x>0 \esp /
\esp \forall y,z\in B(x,\rho_x) \quad
||\tilde{f}_i(y)-\tilde{f}_i(z)||
\leq A_x||y-z||$$
From compacity, we can substract a finite number
of balls which cover $\bar{B}(0,2R)$~:
$$\exists p\in \N \esp / \esp \bar{B}(0,2R) \subset 
\overset{p}{\underset{j=1}{\bigcup}} B(x_j,\rho_{x_j})$$
Let $A=\ds{\sup_i}\, A_{x_i}$ and $\rho=\inv{2}
\min (\f{R}{2},\ds{\min_i}\, \rho_{x_i})$. Let us prove that
$\tilde{f}_i$ is $A$-lipschitzian~:

\no Let $y,z\in \R^n$.
\begin{enumerate}
\item if $||y-z|| \leq \rho$~:
\begin{itemize}
\item if $y,z\in \bar{B}(0,2R)$ then there exists $j\in\{
1,\ldots, p\}$ such that $y,z\in B(x_i,\rho_{x_i})$, and the
conclusion holds.
\item if $y,z\notin \bar{B}(0,R)$ then
$\tilde{f}_i(y)=\tilde{f}_i(z)=0$, and the inequality is still true.
\end{itemize}
All other cases are impossible because $||y-z|| \leq \rho$.
\item if $||y-z|| > \rho$~: \\
Let $M=\ds{\sup_{y,z\in K}} ||\tilde{f}_i(y)-\tilde{f}_i(z)||=
\ds{\sup_{y,z\in \R^n}} ||\tilde{f}_i(y)-\tilde{f}_i(z)||$.
Then~:
$$||\tilde{f}_i(y)-\tilde{f}_i(z)||\leq M \leq
\f{M}{\rho}||y-z||$$
and the conclusion holds if moreover $A$ is choosen larger than
$\f{M}{\rho}$.
\end{enumerate}
\end{proof}
For all $t\in [0,T]$ we have~:
\begin{displaymath}
\begin{split}
||\tilde{x}_u(t)-\tilde{x}_{v}(t)||&=||\, \int_0^t
(\tilde{f}_0(\tilde{x}_u(\tau))-\tilde{f}_0(\tilde{x}_{v}(\tau))
)d\tau \\ & \qquad +\int_0^t \ds{\sum_{i=1}^m} v_{i}(\tau)(
\tilde{f}_i(\tilde{x}_u(\tau))-\tilde{f}_i(\tilde{x}_{v}(\tau))
)d\tau \\ & \qquad 
-\int_0^t \ds{\sum_{i=1}^m} (v_{i}(\tau)-u_i(\tau))
\tilde{f}_i(\tilde{x}_u(\tau)) d\tau \, || \\
&\leq A\int_0^t (1+\ds{\sum_{i=1}^m}|v_{i}(\tau)|)
||\tilde{x}_u(\tau) - \tilde{x}_{v}(\tau)|| d\tau + h_v(t)
\end{split}
\end{displaymath}
where $h_v(t)=\ds{
|| \, \int_0^t \sum_{i=1}^m (v_{i}(\tau)-u_i(\tau))
\tilde{f}_i(\tilde{x}_u(\tau)) d\tau \, ||}$.

\no Set $\ds{M'=\max_i\sup_{x\in\R^n}||\tilde{f}_i(x)||}$.
We get from the \it{Cauchy-Schwarz inequality}~:
$$\forall t\in[0,T]\quad h_v(t)\leq M'\sqrt{T}\, ||v-u||_{L^2}$$
Hence for all $\varepsilon >0$ there exists a neighborhood $U$ of
$u$ in $L^2$ such that~:
$$\forall v\in U\quad \forall t\in[0,T]\quad h_v(t)\leq
\varepsilon $$
Therefore~:
$$\forall t\in [0,T] \quad ||\tilde{x}_u(t)-\tilde{x}_{v}(t)||
\leq A\int_0^t|(1+\sum_{i=1}^m v_i(\tau))|\,
||\tilde{x}_u(\tau)-\tilde{x}_{v}(\tau)||d\tau
+\varepsilon $$
We get from the \it{Gronwall Lemma}~:
$$\ds{\forall t\in [0,T] \quad ||\tilde{x}_u(t)-\tilde{x}_{v}(t)||
\leq \varepsilon \rm{e}^{ A\int_0^t|(1+\sum v_i(\tau))|d\tau }
\leq \varepsilon \rm{e}^{AT+AK\sqrt{T}}}$$
which proves that $(\tilde{x}_{v})$ is
uniformly close to $\tilde{x}_u=x_u$. In particular if the
neighborhood $U$ is small
enough then~: $\forall t\in [0,T] \quad x_{v}(t)\in V$, and hence
$\tilde{x}_{v}=x_{v}$, which ends the proof.
\end{proof}


\subsection{Continuity}
If $v$ and $v_n,n\in \N$ are elements of $L^2([0,T])$, we denote
by $v_n \rightharpoonup v$ the weak convergence of the sequence
$(v_n)$ to $v$ in $L^2$.

\begin{prop} \label{continuity}
Let $u=(u_1,\ldots,u_m)
\in \Omega$ and $x_u$ be the solution of the affine
control system~:
\begin{equation*}
\dot{x}_u=f_0(x_u)+\displaystyle{\sum_{i=1}^{m}u_if_i(x_u)},\esp
x_u(0)=0
\end{equation*}
Let $(u_n)_{n\in \N}$ be a sequence in $L^2([0,T],\R^m)$.
If $u_n \overset{L^2}{\rightharpoonup} u$ then for $n$ large
enough $x_{u_n}$ is well-defined on $[0,T]$ and moreover
$x_{u_n} \tto x_u$ uniformly on $[0,T]$.
\end{prop}

\begin{proof}
The outline of the proof is the same as in Proposition
\ref{existence}.
Let $V$ be a bounded open subset of $\R^n$ such that~:
$\forall t \in [0,T] \quad x_u(t) \in V$. Let $\theta \in
C^\infty (\R^n , [0,1])$ with compact support $K$ such that $\theta
=1$ on $V$. We can assume that $K=\bar{B}(0,R)=\{ x\in \R^n \, /
\, ||x|| \leq R \}$. For $i=0 \ldots m$
we set $\tilde{f}_i=\theta f_i$.
Then it is clear that $x_u$ is also solution of~:
$\dot{x}=\tilde{f}_0(x)+\sum u_i \tilde{f}_i(x)$.
For all $n\in \N$ let $\tilde{x}_{u_n}$ be the solution of
$\dot{\tilde{x}}_{u_n}=\tilde{f}_0(\tilde{x}_{u_n}) +
\displaystyle{\sum_{i=1}^m u_{n,i} \tilde{f}_i
(\tilde{x}_{u_n})}$, $\tilde{x}_{u_n}(0)=0$.
We will prove that if $n$ is large enough then $\tilde{x}_{u_n} =
x_{u_n}$.

For all $t\in [0,T]$ we have~:
\begin{displaymath}
\begin{split}
||\tilde{x}_u(t)-\tilde{x}_{u_n}(t)||&=||\, \int_0^t
(\tilde{f}_0(\tilde{x}_u(\tau))-\tilde{f}_0(\tilde{x}_{u_n}(\tau))
)d\tau \\ & \qquad +\int_0^t \ds{\sum_{i=1}^m} u_{n,i}(\tau)(
\tilde{f}_i(\tilde{x}_u(\tau))-\tilde{f}_i(\tilde{x}_{u_n}(\tau))
)d\tau \\ & \qquad 
-\int_0^t \ds{\sum_{i=1}^m} (u_{n,i}(\tau)-u_i(\tau))
\tilde{f}_i(\tilde{x}_u(\tau)) d\tau \, || \\
&\leq A\int_0^t (1+\ds{\sum_{i=1}^m}|u_{n,i}(\tau)|)
||\tilde{x}_u(\tau) - \tilde{x}_{u_n}(\tau)|| d\tau + h_n(t)
\end{split}
\end{displaymath}
where $h_n(t)=\ds{
|| \, \int_0^t \sum_{i=1}^m (u_{n,i}(\tau)-u_i(\tau))
\tilde{f}_i(\tilde{x}_u(\tau)) d\tau \, ||}$.
The aim is to make $h_n$ uniformly small in $t$, and then to
conclude we use the \it{Gronwall inequality}.

From the hypothesis : $u_n \rightharpoonup u$, we deduce~:
$\forall t\in [0,T] \quad h_n(t) 
\underset{n\to +\infty}{\tto} 0$. Let us prove that $h_n$
tends uniformly to $0$ as $n$ tends to infinity. We need the
following Lemma~:

\begin{lem} \label{dini}
Let $a,b\in \R$ and $E$ be a normed vector space. For all $n\in
\N$ let $f_n:[a,b]\tto E$ be uniformly $\alpha$-hölderian, that is~:
$$\exists \alpha, K>0 \esp / \esp \forall n\in \N \esp \forall
x,y \in [a,b] \quad ||f_n(x)-f_n(y)||\leq K||x-y||^\alpha$$
If the sequence $(f_n)$ converges simply to an application $f$, then
it tends uniformly to $f$.
\end{lem}

\begin{proof}[Proof of the Lemma]
Taking the limit as $n\tto \infty$, it is first clear that $f$
is $\alpha$-hölderian.

Let $\varepsilon >0$ and $a=x_0<x_1<\cdots <x_p=b$ be a
subdivision such that $\forall i \quad
x_{i+1}-x_i<\f{\varepsilon ^\inv{\alpha}}{2K}$.
For all $i$, $f_n(x_i)$ tends to $f(x_i)$, hence~:
$$\exists N\in \N \esp / \esp \forall n\geq N \esp \forall i\in
\{0,\ldots,p \} \quad ||f_n(x_i)-f(x_i)||<\f{\varepsilon }{3}$$
Let $x\in [a,b]$. Then there exists $i$ such that $x\in
[x_i,x_{i+1}]$. Hence~:
\begin{equation*}
\begin{split}
||f_n(x)-f(x)||&\leq 
||f_n(x)-f_n(x_i)||+||f_n(x_i)-f(x_i)||+||f(x_i)-f(x)|| \\
&\leq K||x-x_i||^\alpha +\f{\varepsilon }{3}+K||x-x_i||^\alpha \\
&\leq \varepsilon 
\end{split}
\end{equation*}
\end{proof}
Set~: $M'=\ds{\max_{i} \sup_{x\in \R^n}}||\tilde{f}_i(x)||$. We
get~:
\begin{equation*}
|h_n(x)-h_n(y)|\leq M' \l|\int_y^x (\sum_i |u_{n,i}(\tau)| +
\sum_i |u_i(\tau)|)d\tau \r| 
\end{equation*}
Moreover we get from the \it{Cauchy-Schwarz inequality}~:
\mbox{$\ds{\int_y^x|u| \leq ||u||_{L^2} |x-y|^\inv{2}}$.}
Furthermore the sequence $(u_n)$ converges weakly, hence
is bounded in $L^2$. Therefore there exists a
constant $K$ such that for all $n\in \N$~:
$$|h_n(x)-h_n(y)|\leq K|x-y|^\inv{2}$$
Hence from Lemma \ref{dini} we conclude that the sequence $(h_n)$
tends uniformly to $0$, that is~:
$$\forall \varepsilon >0 \esp \exists N\in \N \esp / \esp
\forall n\geq N \esp \forall t\in [0,T] \quad |h_n(t)|\leq
\varepsilon $$
And hence, if $n\geq N$~:
$$\forall t\in [0,T] \quad ||\tilde{x}_u(t)-\tilde{x}_{u_n}(t)||
\leq A\int_0^t|1+\sum_{i=1}^m v_i(\tau)|\,
||\tilde{x}_u(\tau)-\tilde{x}_{u_n}(\tau)||d\tau
+\varepsilon $$
We get from the \it{Gronwall Lemma}~:
$$\forall t\in [0,T] \quad ||\tilde{x}_u(t)-\tilde{x}_{u_n}(t)||
\leq \varepsilon \rm{e}^{AT+AK\sqrt{T}}$$
which proves that the sequence $(\tilde{x}_{u_n})$ tends
uniformly to $\tilde{x}_u=x_u$. In particular if $n$ is large
enough then~: $\forall t\in [0,T] \quad x_{u_n}(t)\in V$ and hence
$\tilde{x}_{u_n}=x_{u_n}$, which ends the proof.
\end{proof}

\begin{rem}
This Proposition can be found in \cite{Sontag},
but the author uses the following argument~:
if $u_n$ tends weakly to
$0$ then $|u_n|$ tends weakly to $0$, which is not true in
general (take
$u_n(t)=\cos nt$). That is the reason why we need Lemma
\ref{dini}. Otherwise the proof is the same as in \cite{Sontag}.
\end{rem}

To check differentiability in next subsection we will need
the following result~:

\begin{prop} \label{lip}
Let $u\in \Omega$ and $x_u$ be the associated
trajectory. Then for any bounded
neighborhood $U$ of $u$ in $\Omega\subset L^2$ there exists a
constant such that for all $v,w\in U$ and for all $t\in [0,T]$
$$||x_v(t)-x_w(t)||\leq C||v-w||_{L^2}$$
\end{prop}

\begin{proof}
Writing~:
\begin{equation*}
\begin{split}
\dot{x}_v&=f_0(x_v)+\sum_{i=1}^mv_if_i(x_v) \\
\dot{x}_w&=f_0(x_w)+\sum_{i=1}^mw_if_i(x_w) 
\end{split}
\end{equation*}
we get, for all $t\in [0,T]$~:
\begin{equation*}
\begin{split}
||{x}_v(t)-{x}_w(t)||&= ||\, \int_0^t \big(\,
\sum_i (v_i(s)-w_i(s))f_i(x_v(s))+f_0(x_v(s))-f_0(x_w(s)) \\
& \qquad \qquad \qquad
+ \sum_iw_i(s)(f_i(x_v(s))-f_i(x_w(s))) \, \big) ds \,|| \\
&\leq \esp \sum_i \int_0^t |v_i-w_i|\, ||f_i(x_v)||ds + \int_0^t
||f_0(x_v)-f_0(x_w)||ds \\
&\qquad \qquad \qquad + \sum_i \int_0^t
|w_i|\, ||f_i(x_v)-f_i(x_w)||ds
\end{split}
\end{equation*}
Now, if $v$ and $w$ are in a bounded
neighborhood $U$ of $u$ in $L^2$, then according to
Proposition \ref{continuity}, the trajectories $x_v$ and $x_w$
take their values in a compact $K$ that depends only on $U$. The
vector fields $f_0,f_1,\ldots , f_m$ being smooth, we claim
that there exists a constant $M>0$ such that for all $v,w\in U$
and for all $i$
\begin{equation*}
\begin{split}
||f_i(x_v)||&\leq M \\
||f_i(x_v)-f_i(x_u)||&\leq M||x_v-x_u||
\end{split}
\end{equation*}
Endly without loss of generality we can assume that $U$ is
contained in a ball of radius $R$ centered in $O\in L^2$, so that
$$\forall w\in U \quad ||w||_{L^2}\leq R$$
Hence plugging in the upper inequality, and using the
\it{Cauchy-Schwarz inequality}, we obtain~:
$$\forall t\in [0,T] \quad ||x_v(t)-x_w(t)|| \leq A \int_0^t 
||x_v(s)-x_w(s)||ds + B||v-w||_{L^2}$$
where $A$ and $B$ are non negative constants.
Finally we get from the \it{Gronwall Lemma}~:
$$\forall t\in [0,T] \quad ||x_v(t)-x_w(t)|| \leq
C||v-w||_{L^2}$$
with $C=B\rm{e}^{TA}$, which ends the proof.
\end{proof}


\subsection{Differentiability}
Let $u\in\Omega$ and $x_u$ the corresponding solution of the
affine system (\ref{aff0}). We consider the \it{linearized
system} along $x_u$~:
\begin{equation}
\dot{y}_v=A_uy_v+B_uv,\esp
y_v(0)=0,\esp v\in L^2
\end{equation}
where $\ds{A_u(t)=df_0(x_u) +
\sum_{i=1}^m u_idf_i(x_u) }$
and $\ds{B_u(t)=(f_1(x_u),\ldots,f_m(x_u)) }$.
Let $M_u$ be the $n\times n$ matrix solution of
\begin{equation} \label{Mu}
M'_u=A_uM_u, \esp
M_u(0)=Id
\end{equation}
We have~:

\begin{prop} \label{differentiability}
The end-point mapping $E:\begin{array}{rcl} 
\Omega & \tto & \R^n \\
u& \longmapsto & x_u(T) \end{array}$ is $L^2$-Fr\'echet
differentiable, and we have~:
$$ \forall v\in \Omega \quad
dE(u).v=\int_0^TM_u(T)M_u(s)^{-1}B_u(s)v(s)ds$$
\end{prop}

\begin{proof}
Let $u\in L^2([0,T],\R^m)$. Let us prove that $E$ is
differentiable at $u$. Consider a neighborhood $U$ of $0$ in
$\Omega$, and let $v\in U$. Without loss of generality we can assume
that there exists $R>0$ such that for all $v\in U$~: $||v||_{L^2}
\leq R$.
Let $x_u$ (resp. $x_{u+v}$) the solution of the affine system
(\ref{aff0}) with the control $u$ (resp. with the control
$u+v$)~:
\begin{eqnarray}
\dot{x}_{u+v}&=&f_0(x_{u+v})+\sum_{i=1}^m (u_i+v_i)f_i(x_{u+v})
\\
\dot{x}_{u}&=&f_0(x_{u})+\sum_{i=1}^m u_if_i(x_{u})
\end{eqnarray}
We get
$$\dot{x}_{u+v}-\dot{x}_{u}=\sum_{i=1}^m v_if_i(x_{u+v})+
f_0(x_{u+v})-f_0(x_u) +\sum_{i=1}^m u_i(f_i(x_{u+v})-f_i(x_u))$$
Moreover, for all $i=0\ldots m$~:
\begin{equation*}
\begin{split}
f_i(x_{u+v})-f_0(x_u)&=df_i(x_u).(x_{u+v}-x_u)\esp + \\ &\int_0^1
(1-t)d^2f_i(tx_u+(1-t)x_{u+v}).(x_{u+v}-x_u,x_{u+v}-x_u)dt
\end{split}
\end{equation*}
Hence we obtain
\begin{equation} \label{eqdiff}
\dot{\delta}=A_u\delta +B_u\delta +\gamma
\end{equation}
where
$$\delta (t)=x_{u+v}(t)-x_u(t)$$
and
\begin{equation*}
\begin{split}
\gamma (t)=& \sum_{i=1}^m v_i(t) \int_0^1 
df_i(sx_u+(1-s)x_{u+v}).(x_{u+v}-x_u)ds \\ &+ \int_0^1
(1-t)d^2f_0(sx_u+(1-s)x_{u+v}).(x_{u+v}-x_u,x_{u+v}-x_u)ds \\
&+\sum_{i=1}^m u_i(t) \int_0^1
(1-t)d^2f_i(sx_u+(1-s)x_{u+v}).(x_{u+v}-x_u,x_{u+v}-x_u)ds
\end{split}
\end{equation*}
Now for all $v\in U$ we have~: $||v||_{L^2}\leq R$, thus from
Proposition \ref{lip} there exists a compact $K$ in $\R^n$ such
that
$$\forall v\in U \esp \forall s\in [0,1] \quad sx_u(s)
+(1-s)x_{u+v}(s) \in K$$
The $f_i$'s being smooth, we get, using again Proposition
\ref{lip}~:
$$\forall t\in [0,T] \quad ||\gamma(t)||\leq c_1||v||_{L^2}
\sum_{i=1}^m |v_i(t)| + c_2||v||_{L^2}^2(1+\sum_{i=1}^m |u_i(t)|)$$
Now solving equation (\ref{eqdiff}) we obtain
$$\delta(t)=\int_0^t M_u(t)M_u(s)^{-1}B_u(s)v(s)ds
+\int_0^t M_u(t)M_u(s)^{-1}\gamma(s)ds$$
Hence for $t=T$~:
\begin{equation*}
\begin{split}
||x_{u+v}(T)-x_u(T)&-\int_0^T M_u(T)M_u(s)^{-1}B_u(s)v(s)ds|| \\
&\leq c_1||v||_{L^2}\int_0^T\sum_{i=1}^m |v_i(t)|dt
+c_2||v||_{L^2}^2 \int_0^T (1+\sum_{i=1}^m |u_i(t)|)dt \\
&\leq c_3||v||_{L^2}^2
\end{split}
\end{equation*}
Moreover the mapping~:
$\begin{array}{rcl}
L^2 & \tto & \R^n \\
v & \longmapsto & \int_0^T M_u(T)M_u(s)^{-1}B_u(s)v(s)ds
\end{array}$
is linear and continuous. Hence the end-point mapping is Fr\'echet
differentiable at $u$, and its differential at $u$ is this latter
mapping.
\end{proof}

\begin{rem}
Here it was proved that $E$ is differentiable on $L^2$.
It can be found also in \cite{Sontag}. Usually
(see \cite{P}) it is proved to be differentiable on $L^\infty$.
\end{rem}

\begin{rem}
The control $u$ is abnormal and of corank 1 if and only if $\rm{Im
}dE(u)$ is an hyperplane of $\R^n$.
\end{rem}

\begin{prop} \label{contdiff}
With the same assumptions as in Proposition \ref{continuity}, we
have~:
$$u_n \overset{L^2}{\rightharpoonup}
u \Rightarrow dE(u_n) \tto dE(u) \quad \rm{as }n\to +\infty$$
\end{prop}

\begin{proof}
For $s\in [0,T]$, set $N_u(s)=M_u(T)M_u(s)^{-1}$.

\begin{lem}
$N'_u=-N_uA_u, N_u(T)=Id$
\end{lem}

\begin{proof}[Proof of the Lemma]
The matrix $N_uM_u$ is constant as $t$ varies, hence
$(N_uM_u)'=0$. Moreover~: $(N_uM_u)'=N'_uM_u+N_uA_uM_u$, and we
get the Lemma.
\end{proof}

\begin{lem} \label{3.11}
$u_n \overset{L^2}{\rightharpoonup}u \Rightarrow N_{u_n} \tto
N_u$ uniformly on $[0,T]$.
\end{lem}

\begin{proof}[Proof of the Lemma]
For $t\in [0,T]$, we have~:
\begin{equation*}
\begin{split}
N_u(t)-N_{u_n}(t)&=\int_0^t \Big( \, N_{u_n}(s)\big(
df_0(x_{u_n}(s))+\sum_{i=1}^m u_{n,i}(s)df_i(x_{u_n}(s)) \big) \\
& \qquad \qquad
- N_u(s)\big(df_0(x_{u}(s))+\sum_{i=1}^m u_{i}(s)df_i(x_{u}(s)) \big)
\, \Big) \, ds \\
&=\int_0^t \Big( \, \big( N_{u_n}(s)-N_u(s)\big) df_0(x_{u_n}(s))
\\
&\qquad \qquad +
N_u(s)\big( df_0(x_{u_n}(s))-df_0(x_u(s))\big) \\
&\qquad \qquad +\big( N_{u_n}(s)-N_u(s) \big) \sum_{i=1}^m
u_{n,i}(s)df_i(x_{u_n}(s)) \\
&\qquad \qquad +N_u(s)\sum_{i=1}^m
u_{n,i}(s)\big( df_i(x_{u_n}(s))-df_i(x_u(s)) \big) \\
&\qquad \qquad +N_u(s) \sum_{i=1}^m \big(
u_{n,i}(s)-u_i(s) \big) df_i(x_u(s)) \, \Big) \, ds
\end{split}
\end{equation*}
From the hypothesis~: $u_n \rightharpoonup u$, and from
Proposition \ref{continuity}, we get that $x_{u_n}$ tends
uniformly to $x_u$, and hence for all $i$, $df_i(x_{u_n})$ tends
uniformly to $df_i(x_u)$ on $[0,T]$.

\no Secondly, set~: $\ds{h_n(t)=\int_0^t \sum_{i=1}^m \big(
u_{n,i}(s)-u_i(s) \big) N_u(s)df_i(x_u(s))ds }$. Using the same
argument as in the proof of Proposition \ref{continuity}, we
prove that $h_n$ tends uniformly to $0$.

Hence we get the following inequality~:
\begin{equation*}
\begin{split}
\forall \varepsilon >0 \quad \exists N\in \N \esp / \esp &\forall
n\geq N \quad \forall t\in [0,T] \\
& ||N_u(t)-N_{u_n}(t)||\leq
C\int_0^t ||N_u(s)-N_{u_n}(s)||ds + \varepsilon 
\end{split}
\end{equation*}
The \it{Gronwall inequality} gives us~:
$$\forall t\in [0,T] \quad ||N_u(t)-N_{u_n}(t)||\leq \varepsilon
\rm{e}^{CT}$$
and the conclusion holds.
\end{proof}

\begin{lem}
$u_n \overset{L^2}{\rightharpoonup}u \Rightarrow B_{u_n} \tto
B_u$ uniformly on $[0,T]$.
\end{lem}

\begin{proof}[Proof of the Lemma]
From Proposition \ref{continuity}, we know that $x_{u_n}$ tends
uniformly to $x_u$, hence for all $i$, $f_i(x_{u_n})$ tends
uniformly to $f_i(x_u)$, which proves the Lemma.
\end{proof}

We know that the differential of the end-point mapping has the
following form~:
$$\forall v\in L^2([0,T]) \quad dE(u).v=\int_0^T
N_u(s)B_u(s)v(s)ds$$
Therefore from the preceeding Lemmas we get~:
$$\forall v\in L^2([0,T]) \quad dE(u_n).v\tto dE(u).v$$
which ends the proof of the proposition.
\end{proof}


\section{Properties of the value function and of its
level sets}
Let $T>0$ fixed. Consider the affine control system (\ref{aff0})
on $\R^n$ with cost (\ref{cost}).
We denote by $Acc(T)$ the accessibility set in time $T$, that is
the set of points that can be reached from $0$ in time $T$.

\subsection{Existence of optimal trajectories}
The following result is a consequence
of a general result from \cite{LM}, p. 286.

\begin{prop} \label{optimal}
Consider the analytic affine control system in $\R^n$
\begin{equation*}
\dot{x}=f_0(x)+\sum_{i=1}^mu_if_i(x),\esp
x(0)=x_0, \esp x(T)=x_1
\end{equation*}
with cost
$$C(u)=\int_0^T\sum_{i=1}^mu_i^2(t)dt$$
where $T>0$ is fixed and the class $\Omega$ of
admissible controllers is the subset of
$m$-vector functions $u(t)$ in $L^2([0,T],\R^m)$ such that~:
\begin{enumerate}
\item $\forall u\in\Omega\quad$ $x_u$ is well-defined on $[0,T]$.
\item $\exists B_T\esp / \esp\forall u\in \Omega\quad\forall
t\in[0,T]\quad ||x_u(t)||\leq B_T$.
\end{enumerate}
If there exists a control steering $x_0$ to $x_1$, then there
exists an optimal control minimizing the cost steering $x_0$ to
$x_1$.
\end{prop}


\subsection{Definition of the value function}
\begin{defi}
Let $x\in
\R^n$. Define $S:\R^n \tto \R^+ \cup \{+\infty
\}$ by~:
\begin{itemize}
\item If there is no trajectory steering $0$ to $x$ in time $T$,
set ~: $S(x)=+\infty $.
\item Otherwise set~: $S(x)=\inf \{C(u)\esp / \esp u\in
E^{-1}(x) \}$.
\end{itemize}
$S$ is called the value function.
\end{defi}

\begin{defi}
Let $r,T>0$. Define the following level sets~:
\begin{enumerate}
\item $M_r(T)=S^{-1}(r)$.
\item $M_{\leq r}(T)=S^{-1}([0,r])$. 
\end{enumerate}
\end{defi}

\no Combining Proposition \ref{optimal}, arguments of Proposition
\ref{existence} and the fact that the control $u=0$, if
admissible, is minimizing, we get~:

\begin{prop} \label{coroptimal}
Suppose the control $u=0$ is admissible.
Then there exists $r>0$ such that any point of $M_{\leq
r}(T)$ can be reached from $0$ by an optimal trajectory.
\end{prop}

Hence if $r$ is small
enough, $M_r(T)$ (resp. $M_{\leq r}(T)$) is the
set of extremities at time $T$ of minimizing trajectories with
cost equal to $r$ (resp. lower or equal to $r$). It is a
generalization of the (sub)-Riemannian sphere in (sub)-Riemannian
geometry.

\begin{thm} \label{compact}
If $r$ small enough then the subset $M_{\leq r}(T)$ is compact.
\end{thm}

\begin{proof}
First of all, with the same arguments as in
Proposition \ref{existence}, it is easy to see that $M_{\leq
r}(T)$ is bounded if $r$ is small enough.
Now in order to prove that it is closed, consider a
sequence $(x_n)_{n\in \N}$ of points of $M_{\leq r}(T)$
converging to $x\in \R^n$. For each $n$ let $u_n$ be a minimizing
control steering $0$ to $x_n$ in time $T$~: $x_n=E(u_n)$ (the
existence follows from Proposition \ref{coroptimal}). Then
for all $n$, $C(u_n)\leq r$, which means that the sequence
$(u_n)$ is bounded in $L^2([0,T],\R^m)$, and therefore
it admits a weakly converging
subsequence. We can assume that $\ds{u_n
\overset{L^2}{\rightharpoonup }u}$. In particular~: $C(u)\leq r$.
Moreover from Proposition \ref{continuity} we deduce~: $x=E(u)$.
Hence $u$ is a control steering $0$ to $x$ in time $T$ with a
cost lower or equal to $r$. Thus~: $x\in M_{\leq r}(T)$. This
shows that the latter subset is closed.
\end{proof}

\begin{rem}
$M_r(T)$ is not necessarily closed. It is due to the fact that
$S$ can have discontinuities, see Example \ref{exemple}.
\end{rem}


\subsection{Regularity of the value function}
We can now state the main theorem of this section.

\begin{thm} \label{subanalytic}
Consider the analytic affine control system (\ref{aff0})
with cost (\ref{cost}).
Suppose $r$ and $T$ are small enough (so that any trajectory with
cost lower than $r$ is well-defined
on $[0,T]$). Let $K$ be a subanalytic compact subset of
$M_{\leq r}(T)$. Suppose there is no abnormal minimizing geodesic
steering $0$ to any point of $K$. Then $S$ is
continuous and subanalytic on $K$.
\end{thm}

\begin{cor} \label{corsub}
If $r_0$ and $T$ are small enough and if there is no abnormal
minimizer steering $0$ to any point of $M_{\leq r_0}(T)$, then
for any $r$ lower than $r_0$,
$M_r(T)$ and $M_{\leq r}(T)$ are subanalytic subsets of $\R^n$.
\end{cor}

This result generalizes to affine systems a result proved in
\cite{Ag} for sub-Riemannian systems (see also
\cite{AgS},\cite{J}). The main argument
to prove subanalyticity is the same as in \cite{Ag}, i.e.
\it{the compactness of Lagrange multipliers associated to
minimizers}, see Lemma \ref{AV0} below.
\\

If $\Omega=L^2([0,T],\R^m)$, i.e. if trajectories associated to
any control $u$ in $L^2$ are well-defined on $[0,T]$, then any
point of $Acc(T)$ can be joined by a minimizing geodesic.
Theorem \ref{subanalytic} becomes~:

\begin{thm} \label{subanalyticgen}
If $\Omega=L^2([0,T],\R^m)$ and
if there is no abnormal minimizing geodesic, then $S$ is
continuous on $\R^n$ ; moreover $Acc(T)$ is open and
$S$ is subanalytic on any
subanalytic compact subset of $Acc(T)$.
\end{thm}

\begin{proof}[Proof of Corollary \ref{corsub}]
If $r_0$ is small enough then from Theorem \ref{compact}
$M_{\leq r_0}(T)$ is compact. We need a Lemma~:

\begin{lem}
If $r<r_0$ then $M_{\leq
r}(T)$ is contained in the interior of $M_{\leq r_0}(T)$.
\end{lem}

\begin{proof}[Proof of the Lemma]
Let $x$ be a point of $M_{\leq r}(T)$.
From hypothesis, $x$ is the extremity of a
regular geodesic associated to a regular control $u$. Hence $E$
is open in a neighborhood of $u$ in $L^2$. Therefore there exists
a neighborhood $V$ of $x$ such that any point of $V$ can be
reached by trajectories with cost close to $r$~; we can choose
$V$ so that their cost does not exceed $r_0$. Hence $V\subset
M_{\leq r_0}(T)$, which proves that $x$ belongs to the interior
of $M_{\leq r_0}(T)$.
\end{proof}

Let now $K$ be a subanalytic compact subset containing $M_r(T)$
and $M_{\leq r}(T)$. We conclude using Theorem \ref{subanalytic}
and definition of the latter subsets.
\end{proof}

We only prove Theorem \ref{subanalyticgen}. The proof of
Theorem \ref{subanalytic} is similar.

\begin{proof}[Proof of Theorem \ref{subanalyticgen}]
First of all, note that $Acc(T)$ is open. For if $x\in Acc(T)$,
let $u$ be a minimizing control such that $x=E(u)$. From the
assumption, $u$ can not be abnormal. Thus it is normal, and
$dE(u)$ is surjective. Hence from the \it{implicit function
theorem}, $E$ is open in a neighborhood of $u$. Therefore there
exists a neighborhood of $x$ contained in $Acc(T)$, thus
the latter is open.
\\

We first prove the continuity of $S$ on $\R^n$. Take a sequence
$(x_n)$ of points of $\R^n$ converging to $x$. We shall prove
that $S(x_n)$ converges to $S(x)$ by showing that $S(x)$ is the
unique cluster point of the sequence $(S(x_n))$.

\paragraph{First case~:} $x\in Acc(T)$.
Clearly~: $Acc(T)=\underset{r\geq 0}{\cup} M_{\leq r}(T)$, and
moreover~: $r_1<r_2 \Rightarrow M_{\leq r_1}(T) \subset M_{\leq
r_2}(T)$. Hence there exists $r$ such that $x$ and $x_n$ for $n$
great enough are points of $M_{\leq r}(T)$. Now for each $n$
there exists an optimal control $u_n$ steering $0$ to $x_n$, with
a cost $C(u_n)=S(x_n)\leq r$. The sequence $(u_n)$ is bounded in
$L^2$, therefore it admits a weakly converging subsequence. We
can assume that $u_n\rightharpoonup u$. From Proposition
\ref{continuity}, we get~: $x=E(u)$. Let $a$ be a cluster point
of $(S(x_n))_{n\in N}$.
We can suppose that $S(x_n) 
\underset{n\to +\infty}{\tto} a$. From the weak convergence of
$u_n$ towards $u$ we deduce that~: $C(u)\leq a$. Therefore~:
$S(x)\leq a$. Let us prove that actually~: $S(x)=a$. If not, then
there exists a minimizing control $v$ steering $0$ to $x$ with a
cost $b$ strictly lower than $a$. From hypothesis, $v$ is normal,
hence as before $E$ is open in a (strong) neighborhood of $v$ in
$L^2$. It means that points near $x$ can be attained with (not
necessarily minimizing) controls with cost close to
$b$. This contradicts the fact that $S(x_n)$ is close to $a$ if
$n$ is large enough. Hence $a=S(x)$.

\paragraph{Second case~:} $x \notin Acc(T) $.
Then $S(x)=+\infty$.
Let us prove that $S(x_n)\to +\infty$. If not, considering a
subsequence, we can assume that $S(x_n)$
converges to $a$. For each $n$ let $u_n$ be a minimizing control
steering $0$ to $x_n$. Again, the sequence $(u_n)$ is bounded in
$L^2$, hence we can assume that $u_n\rightharpoonup u \in L^2$.
From the continuity of $E$ we deduce~: $x=E(u)$, which is absurd
because $x$ is not reachable. Hence~:
$S(x_n) \underset{n\to +\infty}{\tto} +\infty$.
\\
\\
Let us now prove the subanalyticity property. Let $K$ a compact
subset of $Acc(T)$. Here we use
the first normalization for adjoint vectors (see section 2.2),
that is we choose $p^0=-\inv{2}$ if the extremal is normal.
The following
Lemma asserts that the set of end-points at time $T$ of the
adjoint vectors associated to minimizers steering $0$ to a point
of $K$ is bounded~:
\begin{lem} \label{AV0}
$\{ p_u(T) \esp / \esp E(u)=x_u(T)\in
K, u \rm{ is minimizing} \}$ is a bounded subset of $\R^n$.
\end{lem}

\begin{proof}[Proof of Lemma \ref{AV0}]
If not, there exists a sequence $(x_n)$ of
$K$ such that the associated adjoint vector
verifies~: $||p_n(T)||\underset{n\to +\infty}{\tto} +\infty$.
Substracting a converging subsequence we can suppose that
$x_n \underset{n\to +\infty}{\tto} x$. Now let $u_n$ be a
minimizing control associated to $x_n$, that is~: $x_n=E(u_n)$.
The vector
$p_n(T)$ is a Lagrange multiplier because $u_n$ is minimizing,
hence we have the following equality in $L^2$~:
$$p_n(T).dE(u_n) \overset{L^2}{=} -p^0 u_n $$
Dividing by $||p_n(T)||$ we obtain~:
\begin{equation} \label{lag}
\f{p_n(T)}{||p_n(T)||}.dE(u_n) \overset{L^2}{=}
\f{-p^0}{||p_n(T)||} u_n 
\end{equation}
Actually there exists $r$ such that $K\subset M_{\leq r}(T)$.
Hence~: $C(u_n)\leq r$, and the sequence $(u_n)$ is
bounded in $L^2([0,T],\R^m)$. Therefore
it admits a weakly convergent
subsequence. We can assume that~: $u_n \rightharpoonup u \in
L^2$. Furthermore, the sequence
$\l(\f{p_n(T)}{||p_n(T)||}\r)$ is bounded in $\R^n$, hence up
to a subsequence we have~:
$\f{p_n(T)}{||p_n(T)||} \tto \psi \in \R^n$. Passing to the limit in
(\ref{lag}), and using Proposition \ref{contdiff}, we obtain~:
$$\psi.dE(u)=0 \quad \rm{ where } x=E(u)$$
It means that $u$ is an abnormal control steering $0$ to $x$ in
time $T$. From the assumption, it is not minimizing, hence~:
$C(u)>S(x)$.
On the one part, since $u_n$ is
minimizing, we get from the continuity of $S$ that $C(u_n)\to
S(x)$. On the other part,
from the weak convergence of $(u_n)$ towards $u$ we
deduce that $C(u)\leq S(x)$, and we get a contradiction.
\end{proof}

The previous Lemma asserts that end-points of
adjoint vectors associated
to minimizers reaching $K$ are bounded. We now prove this fact
for initial points of adjoint vectors~:

\begin{lem} \label{AV}
$\{ p_u(0) \esp / \esp E(u)\in
K, u \rm{ is minimizing} \}$ is a bounded subset of $\R^n$.
\end{lem}

\begin{proof}[Proof of Lemma \ref{AV}]
Let $M_u$ be defined as in (\ref{Mu}). From the classical theory
we know that~:
$$p_u(0)=p_u(T)M_u(T)$$
In the same way as in Lemma \ref{3.11} we can prove~:
$$u_n \overset{L^2}{\rightharpoonup }u \Longrightarrow M_{u_n}(T)
\tto M_u(T) \quad \rm{as } n\to +\infty$$
Now if the subset $\{ p_u(0) \esp / \esp E(u)\in
K, u \rm{ is minimizing} \}$ were not bounded, there would exist
a sequence $(u_n)$ such that $||p_{u_n}(0)||\tto +\infty$. Up to
a subsequence we have~: $u_n\rightharpoonup u, x_n=E(u_n)\tto
x\in K$, and with the same arguments as in the previous Lemma, $u$
is minimizing. Then it is clear that
$||p_{u_n}(T)||=||p_{u_n}(0)M^{-1}_{u_n}(T)||
\underset{n\to +\infty}{\tto} +\infty$. This
contradicts Lemma \ref{AV0}.
\end{proof}

Let now $A$ be a subanalytic compact subset of $\R^n$ containing
the bounded subset of Lemma \ref{AV}. Then, if $x\in K$~:
$$S(x)=\inf \{C\rm{o}\phi (p) \esp / \esp p\in (E\rm{o}\phi)^{-1}(x)
\cap A \}$$
(see Definition \ref{fi} for $\Phi$).
Applying Proposition \ref{tamm137} we get
the local subanalyticity of $S$.
\end{proof}

\begin{rem}
In sub-Riemannian geometry (i.e. $f_0=0$) the control $u=0$
steers $0$ to $0$ with a cost equal to $0$, thus is always a
minimizing control. Moreover it is abnormal because
$\rm{Im }dE(0)=\rm{Span }\{f_1(0),\ldots,f_m(0)\}$ has corank
$\geq 1$. Hence hypothesis of Corollary \ref{corsub} \it{is never
satisfied}. That is why the origin must be pointed out. In
\cite{Ag}, Agrachev proves that the sub-Riemannian distance is
subanalytic \it{in a pointed neighborhood of $0$}, and hence that
sub-Riemannian spheres of small radius are subanalytic.

The problem of subanalyticity of the sub-Riemannian distance at
$0$ is not obvious. Agrachev/Sarychev (\cite{AS}) or Jacquet
(\cite{J}) prove this fact under certain assumptions on the
distribution.
In fact for certain dimensions of the state space and
codimensions of the distribution, the absence of abnormal
minimizers (and hence subanalyticity of the spheres) and
non-subanalyticity of the distance  at $0$ are both generic
properties (see \cite{AG}).

Nevertheless for \it{affine systems} with $f_0\neq 0$, the
control $u=0$ (which is always minimizing since $C(u)=0$) is not
in general abnormal. In fact it is not abnormal if and only if
the linearized system along the trajectory of $f_0$ passing
through $0$ is controllable. Such conditions are well-known. For
example we have the following~:
\begin{quote}
If $f_0(0)=0$, set $A=df_0(0), B=(f_1(0),\ldots,f_m(0))$. Then the
control $u=0$ is regular if and only if $rank
(B|AB|\cdots|A^{n-1}B)=n$.
\end{quote}
The regularity property is open, that is~:

\begin{prop}
If $u$ is regular then~:
$$\exists r>0\esp /\esp \forall v\quad ||u-v||_{L^2}
\leq r\Rightarrow v
\rm{ is regular.}$$
\end{prop}

\begin{proof}
If not~: $\forall n\quad \exists v_n\esp /\esp ||u-v_n||_{L^2}
\leq \inv{n}$ and $v_n$ is abnormal. Hence~:
$$\exists p_n, ||p_n||=1\esp /\esp \forall n\quad p_n.dE(v_n)=0$$
Now the sequence $(p_n)$ is bounded in $\R^n$, hence up to a
subsequence $p_n$ converges to $\psi\in \R^n$. On the other
part, $v_n$ converges to $u$ in $L^2$, hence from
Proposition \ref{contdiff} we get~:
$$\psi.dE(u)=0$$
which contradicts the regularity of $u$.
\end{proof}

Hence we can strengthen Corollary \ref{corsub} and state~:

\begin{cor}
Consider the affine system (\ref{aff0}) with cost
(\ref{cost}). If $u=0$ is admissible on $[0,T]$
and is regular, then for any $r$ small enough, $S$ is
continuous on $M_{\leq r}(T)$ and is subanalytic on any
subanalytic compact subset of $M_{\leq r}(T)$. Moreover if $r$ is
small enough then $M_{r}(T)$ and $M_{\leq r}(T)$ are subanalytic
subsets of $\R^n$.
\end{cor}
\end{rem}


\subsection{On the continuity of the value function}
In Theorem \ref{subanalytic} we proved in particular that if
there is no abnormal minimizer then $S$ is continuous on $M_{\leq
r}(T)$. Otherwise it is wrong, as shown in the following example.

\begin{ex}
Consider in $\R^2$ the affine system $\dot{x}=f_0(x)+uf_1(x)$
with
$$f_0=\f{\partial}{\partial x},\esp f_1=\f{\partial}{\partial
y}$$
Fix $T>0$. It is clear that for any $u\in L^2$, $x_u$ is
well-defined on $[0,T]$. We have~:
\begin{equation*}
\begin{split}
x(T)&=T \\
y(T)&=\int_0^T u(t)dt
\end{split}
\end{equation*}
Hence
$$Acc(T)=\{ (T,y) \esp / \esp y\in \R \}$$
The value function takes finite values in $Acc(T)$, and is
infinite outside, thus is not continuous on $\R^n$.
Note that for any control $u$, $dE(u)$ is never
surjective, thus all trajectories are abnormal.
\end{ex}

In the preceding example, $S$ is however continuous in $Acc(T)$.
But this is wrong in general, see the following example.

\begin{ex}[Working example]  \label{exemple}
Consider in $\R^2$ the affine system $\dot{x}=f_0+uf_1$ with
$$f_0=(1+y^2)\f{\partial}{\partial x},\esp
f_1=\f{\partial}{\partial y}$$
Fix $T=1$. The only abnormal trajectory $\gamma$
is associated to $u=0$~: $\gamma(t)=(t,0)$. Let $A=\gamma(1)$ ;
we have $S(A)=0$. The accessibility set at time $1$ is~:
$$Acc(1)=A\cup \{ (x,y)\in \R^2 \esp / \esp x>1 \}$$
Consider now the problem of minimizing the cost
$\ds{C(u)=\int_0^1u^2(t)dt}$. Normal extremals are solutions of~:
$$\begin{array}{rclcrcl}
\dot{x}&=&1+y^2 &\quad& \dot{y}&=&p_y \\
\dot{p}_x&=&0 & & \dot{p}_y&=&-2yp_x
\end{array}$$
Set $p_x=\lambda$.
The area swept by $(x(1),y(1))$ as $\lambda$ varies is represented
on fig. \ref{f3}.

\drawingscale{0.35mm}
\placedrawing{fig3.lp}{$\lambda <0$}{f3}

The level sets $M_r(1)$ of the value function $S$ are represented
on fig. \ref{f5}. The family $(M_r(1))_{r>0}$ is a partition of
$Acc(1)$. Note that the slope of the vector $u_r$ tends to
infinity as $r$ tends to $0$.

\drawingscale{0.5mm}
\placedrawing{fig5.lp}{the level sets of the value function}{f5}

The level sets $M_r(1)$ ramify at $A$, but do not contain this
point, thus they are \it{not closed}. Now we can see that the
value function $S$ \it{is not continuous at $A$, even inside
$Acc(1)$}. Indeed on $M_r(1)$, $S$ is equal to $r$, but at $A$
we have $S(A)=0$.
\\

We can give an equivalent of the value function $S$ near $A$
in the area $(\lambda <0)$ (see fig. \ref{f3}). Computations
lead to the following~:
$$S(x,y)\sim \inv{4}\f{y^4}{x-1}$$

Note that when $y\neq 0$ is fixed, if $x\to 1,x>1$, then
$\lambda\to 
-\infty$. This is a phenomenon of \it{non-properness} due to the
existence of an abnormal minimizer. This fact was already
encountered in sub-Riemannian geometry (see \cite{BT}). In the
next section we explain this phenomenon.
\end{ex}

In this example $A$ is steered from $0$ by the minimizing control
$u=0$. We easily see that the set of minimizing controls steering
$0$ to points near $A$ is not (strongly) compact in $L^2$. In
fact we have the following~:

\begin{thm}  \label{compactcontinu}
Consider the analytic affine control system (\ref{aff0}) with
cost (\ref{cost}). Suppose $r$ and $T$ are small enough. Then $S$
is continuous on $M_{\leq r}(T)$ if and only if the set of
minimizing controls steering $0$ to points of $M_{\leq r}(T)$ is
compact in $L^2$.
\end{thm}

\begin{rem}
In sub-Riemannian geometry the value function $S$ is always
continuous, even though there may exist abnormal minimizers. This
is due to the fact that $S$ is the square of the sub-Riemannian
distance (see for instance \cite{Be}).
Note that in \cite{J} (see also \cite{Ag}) the set of
minimizing controls joining $M_{\leq r}(T)=\overline{B}(0,r)$, $r$
small enough, is proved to be compact in $L^2$.
\end{rem}

\begin{proof}[Proof of Theorem \ref{compactcontinu}]
Suppose $S$ is continuous on $M_{\leq r}(T)$, and let
$(u_n)_{n\in\N}$ be a sequence of minimizing controls steering
$O$ to points $x_n$ of $M_{\leq r}(T)$. From Theorem
\ref{compact} we can assume that $x_n$ converges to $x\in M_{\leq
r}(T)$. Let $u$ be a minimizing control steering $0$ to $x$.
Since $S$ is continuous we get~: $||u_n||_{L^2} \underset{n\to
+\infty}{\tto} ||u||_{L^2}$. The sequence $(u_n)$ is bounded in
$L^2$, hence up to a subsequence it converges weakly to $v\in
L^2$ such that $||v||_{L^2}\leq ||u||_{L^2}$. On the other part
from Proposition \ref{continuity} we get $x=E(v)$. Therefore
$||v||_{L^2}=||u||_{L^2}$ since $u$ is minimizing. Now combining
the weak convergence of $u_n$ to $v$ and the convergence of
$||u_n||_{L^2}$ to $||v||_{L^2}$ we get the (strong) convergence
of $u_n$ to $v$ in $L^2$. This proves the compacity of minimizing
controls since $v$ is minimizing.

The converse is obvious.
\end{proof}


\section{Role of abnormal minimizers}
\subsection{Theorem of tangency}
This analysis is based on the sub-Riemannian Martinet case (see
\cite{BT})~: it was shown that the exponential mapping is not
proper and that in the generic case the sphere is
tangent to the abnormal direction.
This fact is general and we have the following results.

\begin{lem} \label{5.1}
Consider the affine control system (\ref{aff0}) with cost
(\ref{cost}). Assume that there exists a minimizing geodesic
$\gamma$ on $[0,T]$ associated to an unique abnormal minimizing
control $u$ of corank 1,
and that there exists $r>0$ small enough such that
$A=\gamma(T)\in M_{\leq r}(T)$. Denote by $(p_1,0)$ the
projectivized Lagrange multiplier
at $A$. Let $\sigma(\tau)_{0<\tau\leq
1}$ be a curve on $M_{\leq r}(T)$
such that $\ds{\lim_{\tau\to 0}\sigma(\tau)
=A}$. For each $\tau$ we denote by ${\cal P}(\tau)\subset
P(T^*_{\sigma(\tau)}M)$ the set of projectivized
Lagrange multipliers at $\sigma(\tau)$~:
${\cal P}(\tau)=\{(p_u(\tau),p^0_u)\esp /\esp E(u)=\sigma(\tau),
u \rm{ is minimizing} \}$. Then~:
$${\cal P}(\tau)\underset{\tau\to 0}{\tto} \{(p_1,0)\}$$
that is, each Lagrange multiplier of ${\cal P}(\tau)$
tends to $(p_1,0)$ as $\tau\to 0$.
\end{lem}

\begin{proof}
For each $\tau$ let $u_\tau$ a minimizing control
steering $0$ to $\sigma(\tau)$.
For any $\tau\in ]0,1]$ let $(p_\tau(T),p_\tau^0)\in{\cal
P}(\tau)$. Let $(\psi,\psi^0)$ be a cluster point at $\tau=0$~:
there exists a sequence $\tau_n$ converging to $0$
such that $(p_{\tau_n}(T),p^0_{\tau_n})\tto (\psi,\psi^0)$. The
sequence of controls $(u_{\tau_n})$ is bounded in $L^2$, hence up
to a subsequence it converges weakly to a control $v\in L^2$ such
that $C(v)\leq r$. If $r$ is small enough then from Proposition
\ref{coroptimal}, $v$ is admissible. Moreover from Proposition
\ref{continuity} we get~: $E(v)=A$, and the assumption of
the Lemma implies $v=u$. Now write the equality in $L^2$
defining the Lagrange multiplier~:
$$p_{\tau_n}(T).dE(u_{\tau_n})=-p^0_{\tau_n}u_{\tau_n}$$
and passing to the limit we obtain (Proposition \ref{contdiff})~:
$$\psi.dE(u)=-\psi^0u$$
Since $u$ has corank 1, we conclude~:
$(\psi,\psi^0)=(p_1,0)$ in $P(T_A^*M)$.
\end{proof}

Let $\widetilde{E}$ be the end-point mapping for the \it{extended
system} in $\R^n\times \R$~:
\begin{equation} \label{extaff}
\begin{split}
\dot{x}&=f_0(x)+\sum_{i=1}^mu_if_i(x) \\
\dot{x}^0&=\sum_{i=1}^mu_i^2
\end{split}
\end{equation}
If $P\in M_r(T)\subset\R^n$, we denote by $\widetilde{P}=(P,r)$
the corresponding point in the augmented space. In the same way,
we denote by $\widetilde{M}_r(T)$, $\widetilde{M}_{\leq r}(T)$ the
corresponding sets in the augmented space.

\begin{thm} \label{tangent}
Suppose the assumptions of Lemma \ref{5.1} are fulfilled, and set
$r_0=S(A)$. If moreover
$\widetilde{M}_{\leq r}(T)$ is $C^1$-stratifiable near
$\widetilde{A}=(A,r_0)$, then
the strata of $\widetilde{M}_{\leq r}(T)$ are tangent at
$\widetilde{A}$ to the hyperplane $\rm{Im }d
\widetilde{E}(u)$ in $\R^n\times\R$.
If moreover $A\in \overline{M_{r_1}(T)}$, $r_1<r$, then $r_1\geq
r_0$ and the strata of $M_{r_1}(T)$ are tangent at $A$
to the hyperplane $\rm{Im }d{E}(u)$ in $\R^n$, see fig. \ref{f9}.
\end{thm}

\drawingscale{0.5mm}
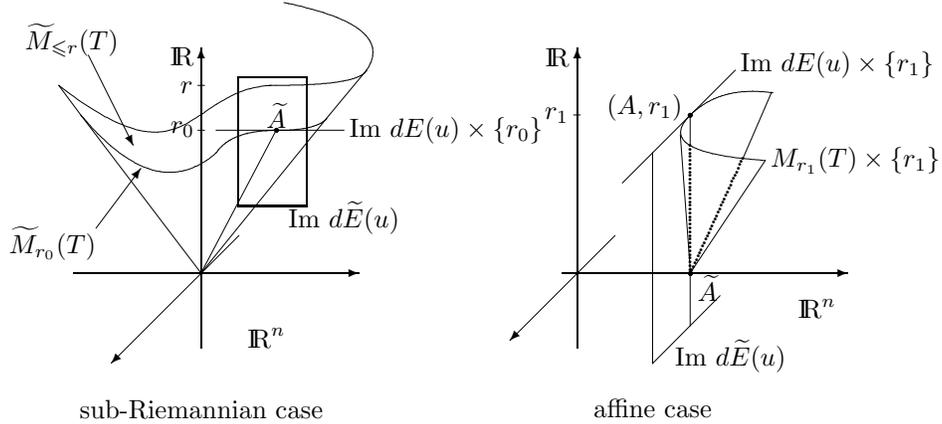
\placedrawing[h]{fig9.lp}{tangency in the augmented space}{f9}

\begin{proof}
Let $N$ be a stratum of $\widetilde{M}_{\leq r}(T)$ of maximal
dimension near $\widetilde{A}$. Let $(\widetilde{\sigma}
(\tau))_{0<\tau\leq 1}$
be a $C^1$ curve on $\overline{N}$ such that $\ds{\lim_{\tau\to
0}\widetilde{\sigma}(\tau)=\widetilde{A}}$, and
$\widetilde{\sigma}(\tau)=(\sigma(\tau),r_\tau)$.
The aim is to prove that $\ds{\lim_{\tau\to 0}
\widetilde{\sigma}'(\tau)\in \rm{Im }d\widetilde{E}(u)}$.
From the assumption on the stratum $\widetilde{N}$,
$\rm{Im }d\widetilde{E}(u_\tau)$ is the tangent space to
$\widetilde{N}$ at $\widetilde{\sigma}(\tau)$.
By definition of the Lagrange multiplier, $(p_\tau(T),p^0_\tau)$
is normal to this subspace.
Moreover $(p_1,0)$ is normal to the hyperplane $\rm{Im
}d\widetilde{E}(u)$. Now from Lemma \ref{5.1} we deduce~:
$\rm{Im }d\widetilde{E}(u_\tau) \underset{\tau \to
0}{\tto} \rm{Im } d\widetilde{E}(u)$. The conclusion is now clear
since $\widetilde{\sigma}'(\tau)\in \rm{Im }d\widetilde{E}(u_\tau)$.

The second part of the Theorem is proved similarly.
\end{proof}

\begin{ex}  \label{martinet}
In \cite{BT} a precise description of the SR sphere in
3-dimensional Martinet case is given. Generically, the abnormal
minimizer has corank 1.
The section of the sphere near the end-point of the abnormal
minimizer with the plane $(y=0)$ is represented on fig. \ref{f67},
(b).

\drawingscale{0.35mm}
\placedrawing[h]{fig67.lp}{}{f67}

In the so-called flat case, the abnormal is not strict, and the
shape of the sphere is represented on fig. \ref{f67}, (a).
In this case, the set of Lagrange multipliers associated to points
near $(-r,0)$
with $z<0$ is bounded. That is why the slope does not converge to
$0$ as $z\to 0, z<0$.

Hence Theorem \ref{tangent} gives a geometric explanation to the
pinching of the generic Martinet sphere near the abnormal
direction.
\end{ex}

\begin{ex}
Consider again the affine system of Example \ref{exemple}. We
proved that the set $M_r(1)$ is tangent at $A$ to the hyperplane
$\rm{Im }dE(u)=\R \f{\partial}{\partial y}$.

Note that, as in the preceding example, computations show that
the branch that ramifies at $A$ is not subanalytic (see fig.
\ref{f5}). In fact, it belongs to the \it{exp-log
category} (see \cite{BLT}).
More precisely this branch has the following graph near $A$~:
$$x=1+F(y,\f{\rm{e}^{-\f{4r}{y^2}}}{y^3})$$
where $F$ is a germ of analytic function at $0$, and we have the
following asymptotic expansion~:
$$x=1+\inv{4r}y^4-3y^2\rm{e}^{-\f{4r}{y^2}}+
\rm{o}(y^2\rm{e}^{-\f{4r}{y^2}})$$
We get the following asymptotics of the value function~:
$$S(x,y)=\inv{4}\f{y^4}{x-1}+\f{y^4}{x-1}\rm{e}^{-\f{y^2}{x-1}}
+\rm{o}(\f{y^4}{x-1}\rm{e}^{-\f{y^2}{x-1}})$$
$S$ is not subanalytic at $A$.
\end{ex}


\subsection{Interaction between abnormal and normal minimizers}
Consider the affine system (\ref{aff0}) with cost (\ref{cost}),
and assume there exists a minimizing geodesic $\gamma$ on $[0,T]$
associated to an unique abnormal control of corank 1. Denote by
$A=\gamma(T)$.

Call \it{normal point} an end-point at time $T$ of a normal
minimizing geodesic. We make the following assumption~:
\begin{quote}
$(H)$ For any neighborhood $V$ of $A$ there exists at least one
normal point contained in $V\cap M_{\leq r}(T)$.
\end{quote}
To describe the normal flow we use the first normalization of
Lagrange multipliers (i.e. $p^0=-\inv{2}$ for normal extremals),
which allows us to define the mapping $\Phi$, see Definition
\ref{fi}. Now set $\rm{exp}=E\rm{o}\Phi$~; it is a generalization
of the (sub)-Riemannian exponential mapping. We have~:

\begin{prop} \label{nonproper}
Under the preceding assumptions the mapping $\rm{exp}$ is not
proper.
\end{prop}

\begin{proof}
Let $(A_n)$ be a sequence of normal points of $M_{\leq r}(T)$
converging to $A$. For each $A_n$ let $(p_n(T),-\inv{2})$ be an
associated Lagrange multiplier. Applying Lemma \ref{5.1} we get~:
$\f{p_n(T)}{\sqrt{||p_n(T)||^2+\inv{4}}} \underset{n\to
+\infty}{\tto} p_1$,
$\f{-\inv{2}}{\sqrt{||p_n(T)||^2+\inv{4}}} \underset{n\to
+\infty}{\tto} 0$,
thus in particular~: $||p_n(T)||\underset{n\to
+\infty}{\tto}+\infty$. Now with the same arguments as in Lemma
\ref{AV} we prove~: $||p_n(0)||\underset{n\to
+\infty}{\tto}+\infty$. By definition~: $A_n=\rm{exp}(p_n(0))$,
hence $\rm{exp}$ is not proper.
\end{proof}

\begin{rem}
Conversely if $\rm{exp}$ is not proper then with the same
arguments as in Lemma \ref{AV0} there exists an abnormal
minimizer. This shows the interaction between abnormal and normal
minimizers. In a sense normal extremals recognize abnormal
extremals. This phenomenon of non-properness is characteristic
for abnormality.
\end{rem}


\subsection{Application~: description of the sub-Riemannian
sphere near an abnormal minimizer for rank 2 distributions}
Let $(M,\Delta,g)$ be a sub-Riemannian structure of rank 2 on an
analytic $n$-dimensional manifold $M$, $n\geq 3$,
with an analytic metric $g$ on $\Delta$.
Our point of view is local and we can assume that $M$ is a
neighborhood of $0\in\R^n$,
and that $\Delta=\rm{Span }\{f_1,f_2\}$ where $f_1,f_2$ are
independant analytic vector fields. Up to reparametrization, the
problem of minimizing the cost (\ref{cost}) at time $T$ fixed is
equivalent to the time-optimal problem with the constraint
$u_1^2+u_2^2\leq 1$.
Let $\hat{\gamma}$ be a reference abnormal trajectory on $[0,r]$,
associated to a control $\hat{u}$ and an adjoint vector
$\hat{p}$. We suppose that $\hat{\gamma}$ is \it{injective}, and
hence without loss of generality we can assume that
$\hat{\gamma}(t)=\rm{exp }tf_1(0)$.

We make the following assumptions~:
\begin{itemize}
\item[$(H_1)$]
Let $K(t)=\rm{Im }dE_t(\hat{u})
=\rm{Span }\{ad^kf_1.f_{2|\hat{\gamma}}, k\geq 0\}$ be the first
Pontriaguine cone along $\hat{\gamma}$. We assume that $K(t)$ has
codimension 1 for any $t\in]0,r]$ and is spanned by the $n-1$
first vectors $ad^kf_1.f_{2|\hat{\gamma}}, k=0\ldots n-2$.
\item[$(H_2)$] $ad^2f_2.f_{1|\hat{\gamma}} \notin K(t)$ along
$\hat{\gamma}$.
\item[$(H_3)$] $f_{1|\hat{\gamma}}\notin
\{ad^kf_1.f_{2|\hat{\gamma}}, k=0\ldots, n-3\}$.
\end{itemize}

Under these assumptions $\hat{\gamma}$ has corank 1. Moreover
from \cite{LS} $\hat{\gamma}$ is minimizing if $r$ is small
enough, and $\hat{u}$ is the unique minimizing abnormal
control steering $0$ to $\hat{\gamma}(r)$. Hence assumptions of
Lemma \ref{5.1} are fulfilled.

Let now $V$ be a neighborhood of $\hat{p}(0)$ such
that all abnormal geodesics starting from $0$ with
$p_\gamma(0)\in V$ satisfy also the assumptions
$(H_1)-(H_2)-(H_3)$.
Note that if $V$ is small enough, they are also injective.
We have, see \cite{AS} and \cite{LS}~:

\begin{prop}
There exists $r>0$ such that the previous abnormal geodesics are
optimal if $t\leq r$.
\end{prop}

\begin{cor}
The end-points of these abnormal minimizers form in the
neighborhood of $\hat{\gamma}(r)$ an
analytic submanifold of dimension $n-3$ if $n\geq 3$, reduced to
a point if $n=3$, contained in the sub-Riemannian sphere
$S(0,r)$.
\end{cor}

Hence in the neighborhood of $\hat{\gamma}(r)$ the
sub-Riemannian sphere $S(0,r)$ splits
into two parts~: the \it{abnormal
part} and the \it{normal part}. To describe $S(0,r)$ near
$\hat{\gamma}(r)$ we have to \it{glue} them together.
If the hypothesis of $C^1$-stratification is fulfilled
then the normal part ramifies
tangently to the abnormal part in the sense of Theorem
\ref{tangent}. This gives us a qualitative description of the
sphere near $\hat{\gamma}(r)$.


\end{document}